\documentclass [11pt,a4paper]{article}

\usepackage{amsmath}
\usepackage{amsthm}
\usepackage{graphics}
\usepackage{latexsym}
\usepackage{amssymb}
\usepackage{enumerate}
\usepackage[latin2]{inputenc}
\usepackage{color}

\newtheorem {theorem}{Theorem}

\newtheorem {lemma}[theorem]{Lemma}

\newtheorem* {theorem*}{Theorem}
\newtheorem* {thm*}{Theorem}
\newtheorem* {lemma*}{Lemma}
\newtheorem* {corollary*}{Corollary}
\newtheorem* {prop*}{Proposition}
\newtheorem* {definition*}{Definition}
\newtheorem* {remark*}{Remark}
\newtheorem* {remarks*}{Remarks}
\def \N {\mathbb N}
\def \Z {\mathbb Z}

\def \R {\mathbb R}

\def \T {\mathbb T}

\def \Ordo {{\cal O}}
\def \ind{1\!\!1}

\newcommand{\abs}[1]{\left|{#1}\right|}

\def\beqs{\begin{eqnarray*}}
\def\eeqs{\end{eqnarray*}}
\def\beq{\begin{eqnarray}}
\def\eeq{\end{eqnarray}}
\def\beas{\begin{eqnarray*}}
\def\eeas{\end{eqnarray*}}
\def\bea{\begin{eqnarray}}
\def\eea{\end{eqnarray}}

\def \prob        {\ensuremath{\mathbf{P}}}
\def \expect      {\ensuremath{\mathbf{E}}}
\def \var         {\ensuremath{\mathbf{Var}}}
\def \cov         {\ensuremath{\mathbf{Cov}}}

\def\pt{\partial_t}
\def\px{\partial_x}

\def\ptt{\pi_{\theta,\tau}}
\def\ttt{\theta,\tau}
\def\txte{\theta(\xi_1+\xi_2)+\tau(\eta_1+\eta_2)}
\def\pxya{\pi(\omega_1)\pi(\omega_2)\pi(\omega_3)}

\def\uo{\underline{\omega}}

\def \OmN {\Omega^{^{_{N}}}}

\def \TN {\T^{^{_{N}}}}
\def \LN {L^{^{_{N}}}}
\def \LNa {L^{*^{_{N}}}}
\def \HN {H^{^{_{N}}}}

\def \XN {X^{^{_{N}}}}

\def \dsum {\displaystyle\sum}

\newcommand{\muN}[1]{\mu^{^{_{N}}}_{#1}}
\newcommand{\nuN}[1]{\nu^{^{_{N}}}_{#1}}
\newcommand{\piN}[1]{\pi^{^{_{N}}}_{{#1}}}
\newcommand{\fN}[1]{f^{^{_{N}}}_{{#1}}}
\newcommand{\hN}[1]{h^{^{_{N}}}_{{#1}}}

\newcommand{\xxi}[1]{{\xi}^{^{_{#1}}}}

\newcommand{\vct}[1]{ \text{\boldmath${#1}$} }
\def \vtt {\vct{\theta}}
\def \vtu {\vct{u}}

\title
{Onsager Relations  and Eulerian Hydrodynamic Limit for
 Systems with Several Conservation Laws}

\author {
{Bálint Tóth \qquad Benedek Valkó}
\\[8pt]
Institute of  Mathematics
\\
Technical University Budapest
}

\begin{document}

\setlength{\baselineskip}{1.23\baselineskip}

\maketitle
\begin{abstract}
We present the derivation of the hydrodynamic limit under Eulerian
scaling for a general class of one-dimensional interacting
particle systems with two or more conservation laws. Following Yau's
relative entropy method it turns out that
in case of more than one conservation laws,
in order that the system exhibit hydrodynamic behaviour,
some particular identities reminiscent of Onsager's reciprocity
relations must hold.
We check validity of these identities for a
wide class of models. It also follows that, as a general rule,
\emph{the equilibrium thermodynamic entropy}
(as function of the densities of the conserved variables)
\emph{is a globally convex Lax entropy}
of the hyperbolic systems of conservation laws arising as hydrodynamic
limit. The Onsager relations arising in this context and its
consequences seem to be novel. As concrete examples we also present a
number  of models modeling deposition (or domain growth) phenomena.
\end{abstract}
\section{Introduction}
\label{section:intro}

We investigate the hydrodynamic behaviour of a very general class of
one dimensional
interacting particle systems with two or more conserved observables.
The systems are not reversible and the hydrodynamic limit under
Eulerian scaling is investigated. We apply Yau's relative entropy
method and obtain validity of the hydrodynamic limit up to the
occurence of the first shock wave in the solution of the limiting
pde. There is no  novelty in the standard steps of the relative
entropy proof, so we only sketch these. The real novelty appears when
it turns out that, in case of more than one conserved quantity,  in
order to complete the relative entropy proof, a class of identities
should hold, relating the macroscopic fluxes appearing in the
hydrodynamic pdes.  These identities are much reminiscent of Onsager's
reciprocity relations. As far as we know these relations have not been
pointed out in the context mathematically rigorous Eulerian
hydrodynamics. We check the validity of these relations for the very
general class of interacting particle systems considered. As a
consequence of the Onsager relations it follows that the systems of
partial differential equations (systems of conservation laws) arising
as hydrodynamic limit are by force of hyperbolic type and the
equilibrium thermodynamic  entropy of the system (as function of the
densities of the conserved quantities) is globally convex Lax entropy
of the hydrodynamic equations. This fact may be not so surprising, as
it is commonly accepted physical fact. So much so that hyperbolic
systems of conservation laws possessing a globally convex Lax entropy
are commonly called \emph{of physical type}, see
\cite{serre}. Nevertheless, as far as we know, this fact has not been
emphasized in the context of mathematically rigorous derivation of
hydrodynamic behaviour.
It is worth noting that given a hyperbolic system of conservation
laws the existence of convex Lax entropies is far from trivial: in the
case of two component systems the local existence of convex Lax
entropies   was established in the very technical work \cite{lax}. In
case of more than two components in general the pdes defining Lax
entropies are overdetermined, so in general Lax entropies do not exist
at all. It turns out from our result that the hyperbolic systems of
conservation laws arising as hydrodynamic limit are of very special
type: they always possess a globally convex Lax entropy, namely the
equilibrium thermodynamic entropy of the system.

Beside the general framework we also present a number of concrete
examples of deposition models with two conserved quantities to which
the general result applies, deriving in this way systems of pdes
(hyperbolic systems of conservation laws) which describe
macroscopically domain growth phenomena in 1+1 dimension.

Our general results are easily extended to more than one dimensions,
however the formalism becomes more complicated. We were mostly
motivated by the (one dimensional) deposition models presented as
concrete examples.

The paper is organized as follows:
In section \ref{section:models} we
present the general formalism and the conditions under which the
hydrodynamic limit is derived.
In section \ref{section:hdl} we state the main results of the paper.
In section \ref{section:examples} we present a number of concrete
examples to which the general framework applies. We hope that the
models introduced in this section could be of interest in the context
of deposition/domain growth phenomena.
In section \ref{section:proof} we sketch the proof of the main result
formulated in section \ref{section:hdl}. The sketchy proof is broken
up into several parts. We only hint at the standard steps of the
relative entropy proof, referring the reader to the original work
\cite{yau} or the monographs \cite{kipnislandim} or \cite{fritz}. The
essential parts of this section are subsections \ref{subs:OR} and
\ref{subs:orconseq} where the Onsager relations and their consequences
are derived. Finally, in section \ref{section:several} we extend the
results formulated in the previous sections from two to arbitrary
number of conserved quantities.


\section{Microscopic models}
\label{section:models}

\subsection{State space, conserved quantities}
\label{subs:statespace}

Throughout this paper we denote by $\T^{^{_{N}}}$ the discrete
tori $\Z/N\Z$, $N\in\N$, and by $\T$ the continuous torus $\R/\Z$.
We will denote the
local
spin state by $S$, we only consider the case
when $S$ is finite. The state space of the interacting particle
system is
\begin{equation*}
\Omega^{^{_{N}}}:=S^{\T^{^{_{N}}}}.
\end{equation*}
Configurations will be denoted
\begin{equation*}
\uo:=(\omega_j)_{j\in \T^{^{_{N}}}}\in\Omega^{^{_{N}}},
\end{equation*}
For sake of simplicity we consider discrete (integer valued) conserved
quantities only.
The two conserved quantities are denoted by
\beqs
\xi: S\rightarrow \Z,
\\
\eta: S\rightarrow \Z,
\eeqs
we also use the notations
$\xi_j=\xi(\omega_j),\, \eta_j=\eta(\omega_j).$
This means that the sums $\sum_j\xi_j$ and $\sum_j\eta_j$ are
conserved by the dynamics.
We assume that the
conserved quantities are different and non-trivial, i.e. the
functions $\xi ,\eta$ and the constant function 1 on $S$ are
linearly independent.


\subsection{Rate function, infinitesimal generator}
\label{subs:rates}

We consider the \emph{rate function}
$r: S\times S\times S\times S\rightarrow \R_+$. The dynamics of the
system consists of elementary jumps effecting nearest neighbour spins,
$(\omega_j,\omega_{j+1})\longrightarrow(\omega'_j,\omega'_{j+1})$,
performed with rate
$r(\omega_j,\omega_{j+1};\omega'_j,\omega'_{j+1})$.

We require that the rate function  $r$
satisfy the following conditions.

\begin{enumerate}[(A)]

\item
\label{cond:cons}
If $r(\omega_1,\omega_2;\omega'_1,\omega'_2)>0$ then
\beq
\label{eq:conserved}
\begin{array}{rcl}
\xi(\omega_1)+\xi(\omega_2)
&=&
\xi(\omega'_1)+\xi(\omega'_2),
\\
\eta(\omega_1)+\eta(\omega_2)
&=&
\eta(\omega'_1)+\eta(\omega'_2).
\end{array}
\eeq

\item \label{cond:irred} For every $K\in [\min \xi,\max \xi], L\in
[\min \eta,\max \eta]$ the set %
\beqs \OmN_{K,L}:= \left\{ \uo\in\OmN:\sum_{j\in\TN}\xi_j=K,
\sum_{j\in\TN}\eta_j=L \right\}\eeqs is an irreducible component
of $\OmN$, i.e. if $\uo, \uo'\in \OmN$ then there exists a series
of elementary jumps with positive rates transforming $\uo$ into
$\uo'$.


\item
\label{cond:staci1}
There exists a probability measure $\pi$ on $S$ such
that for any $\omega_1$, $\omega_2$, $\omega'_1$, $\omega'_2$ $\in S$
\beqs
\pi(\omega_1) \pi(\omega_2) r(\omega_1,\omega_2;\omega'_1,\omega'_2)
=
\pi(\omega_2') \pi(\omega_1') r(\omega_2',\omega_1';\omega_2,\omega_1).
\eeqs

\item
\label{cond:staci2}
We denote
$R(\omega_1,\omega_2):=
\sum_{\omega_1',\omega_2' \in S}
r(\omega_1,\omega_2;\omega'_1,\omega'_2)$.
This is the total jump rate of the nearest neighbour spin pair
$(\omega_1,\omega_2)$. Then for any
$\omega_1,\omega_2,\omega_3 \in S$
\beqs
R(\omega_1,\omega_2)+R(\omega_2,\omega_3)+R(\omega_3,\omega_1)
=
R(\omega_1,\omega_3)+R(\omega_3,\omega_2)+R(\omega_2,\omega_1).
\eeqs

\end{enumerate}

For a precise
formulation of the infinitesimal
generator on $\OmN$ we first define the
map $\Theta_{j}^{\omega',\omega''} :\OmN \rightarrow \OmN$
for every $\omega',\omega''\in S$, $j \in \TN$:
\beqs
\left( \Theta_{j}^{\omega',\omega''} \uo \right)_{i}
=
\left\{
\begin{array}{lcl}
\omega'\quad&\text{ if }& i=j
\\
\omega''\quad&\text{ if }& i=j+1
\\
\omega_i\quad&\text{ if }& i\not=j,j+1.
\end{array}
\right.
\eeqs
The infinitesimal generator of the process defined on $\OmN$
is
\[
\LN f(\uo)=\sum_{j\in \TN} \sum_{\omega',\omega''\in S}
r(\omega_j,\omega_{j+1};\omega',\omega'')
(f( \Theta_{j}^{\omega',\omega''} \uo)-f(\uo)).
\]
We denote by $\XN_t$ the Markov process on the state space $\OmN$ with
infinitesimal generator $\LN$.

\medskip
\noindent
{\bf Remarks:}

\begin{enumerate}[(1)]

\item
Condition (\ref{cond:cons}) implies that $\sum_j\xi_j$ and
$\sum_j\eta_j$
are indeed conserved quantities of the dynamics, while condition
(\ref{cond:irred})
ensures that there are no other hidden conservation laws.

\item
Condition (\ref{cond:irred}) is somewhat imlicit. It seems to us that
it is far not trivial (if not impossible) to formulate explicit
conditions involving the rate functions which would be necessary and
sufficient for irreducibility. However, in the concrete examples
treated in section \ref{section:examples} one can easily check that
irreducibility holds.

\item
Conditions (\ref{cond:cons}), (\ref{cond:irred}) and
(\ref{cond:staci1}) determine the measure $\pi(\omega)$ \emph{up to an
exponential distortion}, that is the probability measures satisfying
these conditions are of the form (\ref{eq:gibbs1}) of the next
subsection.

\item
Conditions (\ref{cond:staci1}) and (\ref{cond:staci2}) imply that the
stationary measures of the  process $\XN_t$  are computable and have
the structure required for hydrodynamic behaviour. See the next
subsection for details. Another consequence of these conditions is
Lemma \ref{lem:OR} which turns out to be of  crucial importance
for hydrodynamic behaviour.

\end{enumerate}


\subsection{Stationary measures, reversed process}
\label{subs:measures}

For every $\theta, \tau \in \R$ let $G(\theta,\tau)$ be the moment
generating function defined below:
\beqs G(\theta,\tau):=\log
\sum_{\omega\in S} e^{\theta \xi(\omega)+\tau \eta(\omega)}
\pi(\omega).
\eeqs
In thermodynamic terms $G(\theta,\tau)$ corresponds to the Gibbs free
energy, see \cite{reichl}.
We define the probability measures
\beq
\label{eq:gibbs1}
\pi_{\theta, \tau}(\omega):=\pi(\omega)
\exp(\theta \xi(\omega)+\tau \eta(\omega)-G(\theta,\tau))
\eeq
on $S$.

Using conditions (\ref{cond:staci1}) and (\ref{cond:staci2}), by very
similar considerations as in
\cite{balazs}, \cite{cocozza}, \cite{rezakhanlou2} or \cite{tothvalko}
one can show that for any $\theta$ and $\tau$ the product
measure
\beqs
\piN{\theta,\tau} :=\prod_{j \in \TN}
\pi_{\theta,\tau}
\eeqs
is stationary for the Markov process on $X^{^{_N}}_t$ on
$\OmN$ with infinitesimal generator $\LN$.
We will refer to these measures as the \emph{canonical}
measures. Since $\sum_j \xi_j$ and $\sum_j \eta_j$ are conserved
the  canonical measures on $\OmN$ are not ergodic.
The conditioned measures defined on $\OmN_{K,L}$ by:
\beqs
\piN{K,L}(\uo):=
\piN{\theta,\tau}\left(\uo
\left|\sum_j \xi_j=K,\, \sum_j \eta_j=L\right.\right)
=
\frac{\piN{\theta,\tau}(\uo)\ind\{\uo\in\OmN_{K,L}\}}
{\piN{\theta,\tau}(\OmN_{K,L})}
\eeqs
are also stationary and due to condtion (\ref{cond:irred}) satisfied
by the rate
functions  they are also ergodic. We shall call these measures the
\emph{microcanonical measures} of our system.
(It is easy to see that the measure $\piN{K,L}$ does
not depend on the values of $\theta, \tau$.)

The elementary movements of the reversed stationary process are
$(\omega_{j-1},\omega_j)$ $\longrightarrow$
$(\omega'_{j-1},\omega'_j)$
with rate $ r(\omega_{j},\omega_{j-1};\omega'_j,\omega'_{j-1})$.
The reversed generator  is
\beqs
\LNa f(\uo)=\sum_{j\in \TN} \sum_{\omega',\omega'' \in S}
r(\omega_j,\omega_{j-1};\omega'',\omega')
(f(\Theta_{j-1}^{\omega',\omega''} \uo)-f(\uo)).
\eeqs
This is the adjoint of the operator $\LN$ with respect to all
microcanonical (and  canonical) measures. I.e.
the reversed process is the same for any $\piN{\theta,\tau}$
or $\piN{K,L}$.


\subsection{Expectations}
\label{subs:expectations}

Expectation, variance, covariance with respect to the
measures  $\piN{\theta,\tau}$
will be denoted by
$\expect_{\theta,\tau}(.)$,
$\var_{\theta, \tau}(.)$,
$\cov_{\theta, \tau}(.)$.

We compute the expectations of the conserved quantities with
respect  to the  canonical measures, as
functions of the parameters $\theta$ and $\tau$:
\beqs
u(\theta,\tau)
&:=&
\expect_{\theta,\tau}(\xi)
=
\sum_{\omega \in S} \xi(\omega) \pi_{\theta,\tau}(\omega)
=
G'_{\theta}(\theta,\tau),
\\
v(\theta,\tau)
&:=&
\expect_{\theta,\tau}(\eta)
=
\sum_{\omega \in S} \eta(\omega) \pi_{\theta,\tau}(\omega)
=
G'_{\tau}(\theta,\tau).
\eeqs
Elementary calculations show, that the matrix-valued function
\beqs
\left(
\begin{array}{cc}
   u'_{\theta} &  u'_{\tau}
\\
   v'_{\theta} &  v '_{\tau}
\\
\end{array}
\right)=\left(
\begin{array}{cc}
   G''_{\theta\theta} &  G''_{\theta\tau}
\\
   G''_{\theta\tau}   &  G''_{\tau\tau}
\\
\end{array}
\right)
=:
G''(\theta,\tau)
\eeqs
is equal to the covariance matrix
$\cov_{\theta, \tau}(\xi,\eta)$ and therefore it is strictly positive
definit. It follows that the function
$(\theta,\tau) \mapsto(u(\theta,\tau),v(\theta,\tau))$
is invertible. We denote the inverse function by
$(u,v)\mapsto(\theta(u,v),\tau(u,v))$.
Actually, denoting by $(u,v)\mapsto S(u,v)$ the convex conjugate
(Legendre
transform) of the strictly convex function $(\theta,\tau)\mapsto
G(\theta,\tau)$:
\beq
\label{eq:thdentropy}
S(u,v):=
\sup_{\theta,\tau}
\big(
u\theta+v\tau-G(\theta,\tau)
\big),
\eeq
we have
\beqs
\theta(u,v)=S'_u(u,v),
\qquad
\tau(u,v)=S'_v(u,v).
\eeqs
In probabilistic terms: $S(u,v)$ is the rate function for joint large
deviations of $(\sum_j\xi_j, \sum_j\eta_j)$. In thermodynamic terms:
$S(u,v)$ corresponds to the equilibrium thermodynamic entropy, see
\cite{reichl}.
Let
\beqs
\left(
\begin{array}{cc}
   \theta'_{u} &  \theta'_{v}
\\
   \tau'_{u} &     \tau'_{v}
\\
\end{array}
\right)=\left(
\begin{array}{cc}
   S''_{uu} &  S''_{uv}
\\
   S''_{uv}   &  S''_{vv}
\\
\end{array}
\right)
=:
S''(u,v).
\eeqs
It is obvious that the matrices $G''(\theta,\tau)$ and $S''(u,v)$ are
strictly positive definit and are inverse of each other:
\beq \label{eq:inverse} G''(\theta,\tau)S''(u,v)=I, \eeq
where either $(\theta,\tau)=(u(\theta,\tau),v(\theta,\tau)) $ or
$(u,v)=(\theta(u,v),\tau(u,v))$. With slight abuse of notation we
shall denote: $\pi_{\theta(u,v),\tau(u,v)}=:\pi_{u,v}$,
$\piN{\theta(u,v),\tau(u,v)}=:\piN{u,v}$,
$\expect_{\theta(u,v),\tau(u,v)}=:\expect_{u,v}$, etc.

We introduce the flux of the conserved quantities. The infinitesimal
generator $\LN$ acts on the conserved quantities as follows:
\beqs
\LN\xi_i=
-
\phi(\omega_{i},\omega_{i+1})
+
\phi(\omega_{i-1},\omega_{i})
=:
-\phi_{i} + \phi_{i-1},
\\
\LN\eta_i=
-
\psi(\omega_{i},\omega_{i+1})
+
\psi(\omega_{i-1},\omega_{i})
=:
-\psi_{i} + \psi_{i-1},
\eeqs
where
\beq
\label{eq:phipsidef}
\begin{array}{rrl}
\phi(\omega_1,\omega_2)
&:=&
\dsum_{\omega'_1,\omega'_2\in S}
r(\omega_1,\omega_2;\omega'_1,\omega'_2)
(\xi(\omega'_2)-\xi(\omega_2))
\\[8pt]
&=&
\dsum_{\omega'_1,\omega'_2\in S}
r(\omega_1,\omega_2;\omega'_1,\omega'_2)
(\xi(\omega_1)-\xi(\omega'_1)),
\\[18pt]
\psi(\omega_1,\omega_2)
&:=&
\dsum_{\omega'_1,\omega'_2\in S}
r(\omega_1,\omega_2;\omega'_1,\omega'_2)
(\eta(\omega'_2)-\eta(\omega_2))
\\[8pt]
&=&
\dsum_{\omega'_1,\omega'_2\in S}
r(\omega_1,\omega_2;\omega'_1,\omega'_2)
(\eta(\omega_1)-\eta(\omega'_1)).
\end{array}
\eeq
We shall denote the
expectations of these functions with respect to the canonical measure
$\pi^{^{_{2}}}_{u,v}$  by
\beq
\label{eq:PhiPsidef}
\begin{array}{rrl}
\Phi(u,v)
&:=&
\expect_{u,v} (\phi)
\\
&=&
\dsum_{\substack{\omega_1,\omega_2,\\\omega'_1,\omega'_2}\in S}
r(\omega_1,\omega_2;\omega'_1,\omega'_2)
(\xi(\omega'_2)-\xi(\omega_2))
\pi_{u,v}(\omega_1) \pi_{u,v}(\omega_2),
\\[25pt]
\Psi(u,v)
&:=&
\expect_{u,v} (\psi)
\\
&=&
\dsum_{\substack{\omega_1,\omega_2,\\\omega'_1,\omega'_2}\in S}
r(\omega_1,\omega_2;\omega'_1,\omega'_2)
(\eta(\omega'_2)-\eta(\omega_2))
\pi_{u,v}(\omega_1) \pi_{u,v}(\omega_2).
\end{array}
\eeq
The first derivative matrix of the fluxes $\Phi$ and $\Psi$ will be
denoted
\beq
\label{eq:fluxderiv}
D(u,v):=
\left(
\begin{array}{cc}
\Phi'_u&\Phi'_v
\\
\Psi'_u&\Psi'_v
\end{array}
\right)
\eeq

As a general convention,  if $\delta:S^m\to\R$ is a local function
then its expectation  with respect to
the canonical measure $\pi^{^{_{m}}}_{u,v}$ is denoted
by
\beqs
\Delta(u,v)
:=
\expect_{u,v} (\delta)
=
\sum_{\omega_1,\dots,\omega_m\in S^m}
\delta(\omega_1,\dots,\omega_m)
\pi_{u,v}(\omega_1) \cdots
\pi_{u,v}(\omega_m).
\eeqs
%


\section{The hydrodynamic limit}
\label{section:hdl}

We will show, applying Yau's relative entropy method, that under
Eulerian scaling the local densities of the conserved quantities $
u(t,x),\, v(t,x)$ evolve according to the following  system of partial
differential equations:
\beq
\label{eq:pde}
\left\{\begin{array}{rcl}
  \pt u+\px \Phi(u,v)&=&0
\\
  \pt v+\px \Psi(u,v)&=&0.
\\
\end{array}
\right.
\eeq
It also turns out from our proof (more precisely as a consequence of
the Onsager relations proved in Lemma \ref{lem:OR}) that the
\emph{systems of conservation laws}
(\ref{eq:pde}) arising as hydrodynamic limits are necessarily
\emph{of hyperbolic type}
and the equilibrium thermodynamic entropy function
$(u,v)\mapsto S(u,v)$  is a (very special) globally
\emph{convex Lax entropy}
for the system (\ref{eq:pde}).
(See \cite{serre} or \cite{smoller} for the pde notions used.)
This  may be not so surprising, as it is commonly
accepted fact and drops out automatically, without any computations in
some particular model systems investigated so far.
Nevertheless, we have not found a general statement or proof of this
fact in the hydrodynamic limit literature.

\subsection{Notations}
\label{subs:notations}

For the proper formulation of our results we need some more
notations. Let $u(t,x), v(t,x), t\in [0,T], x\in S$ be a smooth
solution of  (\ref{eq:pde}) (more precisely: let it be twice
continuously differentiable in both variables). We shall use the
notations
\beqs
\theta(t,x)&:=&\theta(u(t,x), v(t,x))
\\
\tau(t,x)&:=&\tau(u(t,x), v(t,x)).
\eeqs

The \emph{true distribution} of
the Markov process $\XN_s$ at macroscopic time $t$, i.e.,  at
microscopic time $Nt$ is
\begin{equation}
\label{eq:mudef}
\muN{t} := \muN{0} \exp \left\{ N t \LN
\right\}.
\end{equation}
%
%

The true distribution will be compared to the following \emph{time
dependent reference measure} (also called local equilibrium) on $\OmN$:
\beq
\label{eq:nudef}
\nuN{t}:=\prod_{j\in \TN}
\pi_{u(t,\frac{j}{N}),v(t,\frac{j}{N})}.
\eeq
This measure is not stationary (unless $u$ and $v$ are
constant), and the local densities of the conserved quantities are
discrete approximations of the functions $u(t,x),  v(t,x)$.

We shall use a stationary measure $\piN{}:=\piN{0,0}$ on $\OmN$
as an \emph{absolute  reference measure}.
The Radom-Nikodym
derivatives of the true distribution and the time dependent
reference measure with respect to the absolute referencee measure are
denoted as follows:
\beq
\label{eq:RadNikdef}
\nonumber
\hN{t}
&:=&
\frac{d\muN{t}}{d\piN{}}(\uo)=\exp\{N t \LNa\}\hN{0}.
\\
\fN{t}
&:=&
\frac{d\nuN{t}}{d\piN{}}(\uo)
\\
\nonumber
&=&
\prod_{j\in \TN} \exp\{ \xi(\omega_j)
\theta(t,\frac{j}{N})+\eta(\omega_j)
\tau(t,\frac{j}{N})-G(\theta(t,\frac{j}{N}),\tau(t,\frac{j}{N}))\}
\eeq

\subsection{The main result}
\label{subs:result}

Our aim is to prove that if $\muN{0}$ is close to $\nuN{0}$ in the
sense of relative entropy, then $\muN{t}$ stays close to $\nuN{t}$
in the same sense uniformly for $t\in[0,T]$. If we consider two
different pairs of smooth solutions $(u_i(t,x),v_i(t,x)), \,\,\,i=1,2$
of (\ref{eq:pde}) it is a simple exercise to show that the relative
entropy of the two time dependent reference measures is of order
$\asymp N$. This suggests that one should prove
\beq
\label{eq:main}
\HN(t):=H(\muN{t}|\nuN{t})=o(N),
\eeq
uniformly for $t\in[0,T].$

\begin{theorem*}
Consider an interacting particles system model defined as in the
previous section which satisfies conditions
(\ref{cond:cons}), (\ref{cond:irred}), (\ref{cond:staci1}) and
(\ref{cond:staci2}).  Let $\Phi(u,v)$ and $\Psi(u,v)$ be defined as in
(\ref{eq:PhiPsidef}).
\\
(i) The system of conservation laws (\ref{eq:pde}) is hyperbolic
in the domain $(u,v)\in
(\min\xi,\max\xi)\,\times\,(\min\eta,\max\eta)$. Furthermore,  the
equilibrium thermodynamic entropy $(u,v)\mapsto S(u,v)$ is a
globally convex Lax entropy for the system (\ref{eq:pde}).
\\
(ii)
Let $[0,T]\times\T\ni (t,x) \mapsto (u(t,x), v(t,x))$ be a
smooth solution of  (\ref{eq:pde}), and let $\muN{t}$ and  $\nuN{t}$
be the
measures defined in (\ref{eq:mudef}), respectively,
(\ref{eq:nudef}). Then, if
\beqs
H(\muN{0} |\piN{})=\Ordo(N),
\eeqs
and (\ref{eq:main}) holds
for $t=0$ then it will hold uniformly for $t\in[0,T]$.
\end{theorem*}

\noindent
{\bf Remark:}
Part (i) of the Theorem is commonly accepted fact.
In some  particular models investigated it simply drops out
without any computation.
However, we do not know about any explicit formulation (or proof) of
the  \emph{general fact} stated here.

\smallskip

From part (ii) of the Theorem, by applying the entropy inequality in a
standard way (comparing the true measure $\muN{t}$ with the local
equilibrium reference measure $\nuN{t}$) one gets the following
corollary:

\begin{corollary*}
Under the conditions of the Theorem, for any $t \in[0,T]$,  the
following limits  hold  as $N\rightarrow \infty$:
\\
(i) For any smooth test function
$g: \T \rightarrow\R$%
\beqs
\frac1N \sum_{j\in \TN}
g({j}/{N})\xi_j(t){\buildrel\prob\over\longrightarrow}
\int_\T g(x)u(t,x)\,dx,
\\
\frac1N \sum_{j\in \TN}
g({j}/{N})\eta_j(t){\buildrel\prob\over\longrightarrow}
\int_\T g(x)v(t,x)\,dx.
\eeqs
\\
(ii)
The asymptotics of the relative entropy of the true distribution
$\muN{t}$ with respect to the absolute
reference measure $\piN{u_0,v_0}$ is
\beq
\label{eq:entrolim}
\N^{-1}H(\muN{t}\,|\,\piN{u_0,v_0})
\to
\int_\T \big(S(u(t,x),v(t,x))-S(u_0,v_0)\big)\,dx,
\eeq
where $S(u,v)$ is the thermodynamic entropy defined in
(\ref{eq:thdentropy}).

\end{corollary*}

\noindent
{\bf Remark:} 
Note that since $S(u,v)$ is Lax entropy of the pde (\ref{eq:pde})
the right hand side of (\ref{eq:entrolim}) does not change in
time as long as the solution $(u(t,x),v(t,x))$ of  (\ref{eq:pde})
is smooth, and starts to decrease when the first shock
appears. This means that the relative entropy
$H(\muN{t}\,|\,\piN{u_0,v_0})$ decreases by $o(N)$ before
the appearence of the first shock in the system.


\section{Examples: deposition models}
\label{section:examples}

If we fix the size of the spin space, then we have finitely many
equations from the conditions on the rate function, thus we can
get a finite-parameter family  of models. The smallest value of
$\abs{S}$, for which there exists a proper model is 3, since we
need to have two different non-trivial conserved quantities.
We present two concrete examples: one with $\abs{S}=3$, one with
$\abs{S}=4$. A third example with $\abs{S}=\infty$, to which the
Theorem applies with some modifications, is described in
\cite{tothwerner}.

Our concrete examples are \emph{deposition models}. $\eta:S\to\N$,
and  $\xi:S\to\Z$. $\eta_j$, respectively, $\xi_j$ are interpreted as
particle occupation number, respectively, (negative) discrete gradient
of deposition height. The dynamical driving mechanism is such that
\\
(i)
The deposition height growth is influenced by the local particle
density. Typically: growth is enhanced by higher particle densities.
\\
(ii)
The particle motion is itself influenced by the deposition
profile. Typically: particles are pushed in the direction of the
negative gradient of the deposition height.

It is natural to assume left-right symmetry of the models. This is
realized in the following way. There is an involution
\[
R:S\to S, \qquad
R\circ R = Id
\]
which acts on the conserved quantities and the jump rates as follows:
\beq
\nonumber
&
\eta(R\omega)=\eta(\omega),
\qquad
\xi(R\omega)=-\xi(\omega),
\\
\label{eq:refl}
&
r(R\omega_2,R\omega_1;R\omega'_2,R\omega'_1)
=
r(\omega_1,\omega_2;\omega'_1,\omega'_2).
\eeq

Correspondingly on the macroscopic level we shall use the
traditional notation $\rho(t,x)$ (instead of $v(t,x)$ of the general
formulation) and $u(t,x)$. The limiting partial differential equations
will be also invariant under the left-right reflection symmetry:
$\big(\rho(t,x),u(t,x)\big) \mapsto \big(\rho(t,-x),-u(t,-x)\big)$.

\subsection{A model with $\abs{S}=3$}
\label{subs:threestates}

The state space is $S=\{-1,0,1\}$. The left-right reflection symmetry
is implemented by $R:S\to S$, $R\omega=-\omega$. The two conserved
quantities
are $\xi(\omega)=\omega$ (the spin itself) and
$\eta(\omega)=1-\abs{\omega}$ (the number of zeros).
(It is easy to see that up to linear combinations these are
the only two conserved quantities we can define on $S$.)
From condition  (\ref{cond:cons}) it follows that
$r(\omega_1,\omega_2;\omega'_1,\omega'_2)>0$ only if
$\omega'_1=\omega_2$ and
$\omega'_2=\omega_1$.
I.e. the dynamics comsists of exchanges of
nearest neighbour spins. It follows  that, without any restriction on
the rates, condition (\ref{cond:staci1})
is satisfied with any probability measure $\pi$ on $S$. Our natural
parametrization is
\[
\pi_{\rho,u}(0)=\rho, \quad
 \pi_{\rho,u}(\pm1)=\frac{1-\rho\pm u}{2},
\]
with the parameter range
$\{(\rho,u):\rho\in[0,1],\,\,u\in[-1,1],\,\,\rho+|u|\le1\}$.

Condition (\ref{cond:staci2}) is fulfilled if and
only if
\beqs
&&
r(1,-1;-1,1)-r(-1,1;1,-1)
\\
&&
\phantom{MMM}=r(1,0;0,1)-r(0,1;1,0)+r(0,-1;-1,0)-r(-1,0;0,-1)
\eeqs
holds. The reflection symmetry condition (\ref{eq:refl}) reads
\[
r(1,0;0,1)=r(0,-1;-1,0),
\qquad
r(0,1;1,0)=r(-1,0;0,-1).
\]
These conditions leave us with
\[
\begin{array}{rcrrcl}
 r(1,-1;-1,1)&=&a, & \qquad r(-1,1;1,-1)&=&2c+a,
\\
 r(0,-1;-1,0)&=&b, &\qquad r(-1,0;0,-1)&=&c+b,
\\
 r(1,0,0,1)&=&b, & \qquad r(0,1,1,0)&=&c+b,
\\
\end{array}
\]
where $a,b \ge0$ and $c\ge\max\{-b,-a/2\}$ are free parameters.
Without loss of generality we may choose $c\ge0$ (otherwise,
rename $\tilde\omega:=-\omega$). It is easy to check that ondition
(\ref{cond:irred}) is satisfied iff $(a+2c)(b+c)>0$. We are not
interested in the $c=0$ case, since that defines the reversible
process which would imply diffusive rather than hyperbolic
(Eulerian) scaling. By fixing an appropriate time scale we choose
$c=1$. It is easy to compute the microscopic fluxes $\phi_j$ and
$\psi_j$ given by formula (\ref{eq:phipsidef}):
\beqs \phi_j &=& \frac12\left(\omega_j - \omega_{j+1}\right)
\left((\omega_j - 1)(1 + \omega_{j+1}) -2 a \omega_j \omega_{j+1}
+ 2 b \left(1 + \omega_j \omega_{j+1}\right)\right),
\\[8pt]
\psi_j
&=&
b(\omega^2_{j+1}-\omega^2_j)
+\frac12 (1-\omega_j)(1+\omega_{j+1})(\omega_j+\omega_{j+1})
\eeqs
The macrosscopic fluxes are computed with formula
(\ref{eq:PhiPsidef}.  Inserted in (\ref{eq:pde}) this leads
to the hydrodynamic equation:
\beq
\label{eq:leroux}
\left\{
\begin{array}{lcl}
  \pt \rho + \px ( \rho u)     &=&0
\\
  \pt u    + \px ( \rho + u^2) &=&0.
\end{array}
\right.
\eeq
This system is known in the pde community as \emph{Leroux's equation}.
The system has some very special features: it belongs to
the so-called Temple class and it was much investigated. For details
see \cite{serre}. Validity of this pde in the hydrodynamic limit
\emph{up to the occurence of shocks} follows from our general
Theorem.

\smallskip
\noindent
{\bf Remark:}
It is an easy exercise to see that all models with $\abs{S}=3$
satisfying the general conditions (\ref{cond:cons},
\ref{cond:irred}, \ref{cond:staci1}, \ref{cond:staci2}), without
the extra assumption of left-right reflection symmetry,
are essentially equivalent, in the sense that
in the hydrodynamic limit they lead to pde-s which can be transformed
to (\ref{eq:leroux}) by linear combinations of the functions involved.

\subsection{A finite bricklayer model}
\label{subs:fourstates}

In the following example we give a finite version of the infinite
bricklayers model introduced in \cite{tothwerner}.
Let $S=\{0,1\}\times\{-1,1\}$. The elements of S will be denoted
$\omega=:(n,z)$. Left-right reflection symmetry is implemented as
$R:S\to S$, $R(n,z)=(n,-z)$. The conserved quantities are
$\xi(\omega)=z$
and
$\eta(\omega)=n$. Condition (\ref{cond:cons}) leaves twenty (possibly)
non-zero rates. Due to the left-right reflection symmetry conditions
eight pairs of rates are equal.  Using the notation
\[
r(\omega_1,\omega_2;\omega'_1,\omega'_2) = r\left(
\substack{n_1^{\phantom{a}}\\z_1^{\phantom{a}}} ,
\substack{n_2^{\phantom{a}}\\z_2^{\phantom{a}}} ;
\substack{n^\prime_1\\z'_1} , \substack{n'_2\\z'_2}, \right)
\]
in the following table we list
the (possibly) non-zero rates, parametrized by twelve nonnegative
parameters.
\beqs
&
r\left(
\substack{0\\-}
,
\substack{0\\+}
;
\substack{0\\+}
,
\substack{0\\-}
\right)
=
a,
\qquad
r\left(
\substack{0\\+}
,
\substack{0\\-}
;
\substack{0\\-}
,
\substack{0\\+}
\right)
=
b,
\\
&
r\left(
\substack{1\\-}
,
\substack{1\\+}
;
\substack{1\\+}
,
\substack{1\\-}
\right)
=
c,
\qquad
r\left(
\substack{1\\+}
,
\substack{1\\-}
;
\substack{1\\-}
,
\substack{1\\+}
\right)
=
d,
\\
&
r\left(
\substack{0\\-}
,
\substack{1\\+}
;
\substack{0\\+}
,
\substack{1\\-}
\right)
=
e,
\qquad
r\left(
\substack{1\\-}
,
\substack{0\\+}
;
\substack{1\\+}
,
\substack{0\\-}
\right)
=
e,
\\
&
r\left(
\substack{0\\+}
,
\substack{1\\-}
;
\substack{0\\-}
,
\substack{1\\+}
\right)
=
f,
\qquad
r\left(
\substack{1\\+}
,
\substack{0\\-}
;
\substack{1\\-}
,
\substack{0\\+}
\right)
=
f,
\\
&
r\left(
\substack{0\\-}
,
\substack{1\\-}
;
\substack{1\\-}
,
\substack{0\\-}
\right)
=
p,
\qquad
r\left(
\substack{1\\+}
,
\substack{0\\+}
;
\substack{0\\+}
,
\substack{1\\+}
\right)
=
p,
\\
&
r\left(
\substack{0\\+}
,
\substack{1\\+}
;
\substack{1\\+}
,
\substack{0\\+}
\right)
=
q,
\qquad
r\left(
\substack{1\\-}
,
\substack{0\\-}
;
\substack{0\\-}
,
\substack{1\\-}
\right)
=
q,
\\
&
r\left(
\substack{0\\+}
,
\substack{1\\-}
;
\substack{1\\+}
,
\substack{0\\-}
\right)
=
r,
\qquad
r\left(
\substack{1\\+}
,
\substack{0\\-}
;
\substack{0\\+}
,
\substack{1\\-}
\right)
=
r,
\\
&
r\left(
\substack{0\\-}
,
\substack{1\\+}
;
\substack{1\\-}
,
\substack{0\\+}
\right)
=
s,
\qquad
r\left(
\substack{1\\-}
,
\substack{0\\+}
;
\substack{0\\-}
,
\substack{1\\+}
\right)
=
s,
\\
&
r\left(
\substack{0\\-}
,
\substack{1\\+}
;
\substack{1\\+}
,
\substack{0\\-}
\right)
=
x,
\qquad
r\left(
\substack{1\\-}
,
\substack{0\\+}
;
\substack{0\\+}
,
\substack{1\\-}
\right)
=
x,
\\
&
r\left(
\substack{0\\+}
,
\substack{1\\-}
;
\substack{1\\-}
,
\substack{0\\+}
\right)
=
y,
\qquad
r\left(
\substack{1\\+}
,
\substack{0\\-}
;
\substack{0\\-}
,
\substack{1\\+}
\right)
=
y.
\eeqs
All the other jump rates are zero.

It is easy to check that condition (\ref{cond:staci1}) imposes
\beq
\label{eq:r=s}
r=s
\eeq
and no other restriction. It also follows that the measures
$\pi_{\rho,u}$ are of the product form
\beq
\label{eq:pinz}
\pi_{\rho,u}(n,z)=
\big(n\rho+(1-n)(1-\rho)\big)\frac{1+z u}{2},
\quad
n=0,1,\,\,\,\,\,z=+,-,
\eeq
with the parameters $\rho\in(0,1)$, $u\in(-1,+1)$.

Another straightforward computation shows that condition
(\ref{cond:staci2}) reads
\beq
\label{eq:turo}
\begin{array}{c}
c+f+p+y = d+e+q+x
\\
a+f+q+y = b+e+p+x
\end{array}
\eeq
So, we we are left with  a nine-parameter family of models.
Given the  formulas (\ref{eq:phipsidef}) we compute the fluxes
$\phi_j$ and $\psi_j$. Using the conditions (\ref{eq:r=s}) and
(\ref{eq:turo}) eventually we get
\beqs
2\phi_j
&=&
\phantom{+}
(b-a)
+
(p-q)(n_j+n_{j+1})
-
(b+a)(z_{j+1}-z_j)
\\
&&
+
(a+b-e-f-x-y)(n_j+n_{j+1})(z_{j+1}-z_j)
+
(a-b)z_{j+1}z_j
\\
&&
-
(a+b+c+d-2e-2f-2x-2y) n_jn_{j+1}(z_{j+1}-z_j)
\\
&&
-
(p-q)
(n_j+n_{j+1})z_{j+1}z_j
\\
4\psi_j
&=&
-
(p+q+r+s+x+y)(n_{j+1}-n_j)
\\
&&
+
(p-q)(n_j+n_{j+1})(z_j+z_{j+1})
+
(y-x)(n_{j+1}-n_j)(z_{j+1}-z_j)
\\
&&
-2(p-q)n_jn_{j+1}(z_j+z_{j+1})
\\
&&
-(p+q-r-s-x-y) (n_{j+1}-n_j)z_{j+1}z_j
\eeqs
The macroscopic fluxes are again explicitly computable.
From (\ref{eq:PhiPsidef}) and (\ref{eq:pinz}) we
get
\beqs
\Phi(\rho,u)
&=&
\big((p-q)\rho-(a-b)/2\big)\big(1-u^2\big),
\\
\Psi(\rho,u)
&=&
(p-q)\rho(1-\rho)u.
\eeqs
Without loss of generality we may assume $p-q\ge0$ (otherwise rename
the microscopic variables $\tilde n_j:=n_j$, $\tilde
z_j:=-z_j$). Further on, $p=q$ leads to diffusive rather then
hyperbolic (Eulerian) scaling of the particle density, so we are
interested in the $p>q$ cases only.
By setting the appropriate  time scale we can choose $p-q=1$ and denote
$\gamma:=\frac{a-b}{2(p-q)}$. So, eventually we get the system of pdes
\beq
\label{eq:brickl}
\left\{
\begin{array}{l}
\pt\rho+\px\big(\rho(1-\rho)u\big)=0
\\[8pt]
\pt u  +\px\big((\rho-\gamma)(1-u^2)\big)=0
\end{array}
\right.
\eeq

In \cite{popkovschuetz} another family of four-state models with two
conserved quantities, the so-called \emph{two channel traffic models}
are analyzed.
These models also satisfy  conditions
(\ref{cond:cons}), (\ref{cond:irred}), (\ref{cond:staci1}) and
(\ref{cond:staci2}). As a consequence our general Theorem is applicable
to  the two channel traffic models, too.

About the relation of our four state deposition models (treated in
this subsection) and the  two channel traffic models treated in
\cite{popkovschuetz}:
Due to the different symmetry conditions imposed
--- we impose the left-right reflection symmetry described in the
first paragraph of this section, while in \cite{popkovschuetz}
symmetry between the two traffic channels is imposed ---
the two families of models do not intersect.
The one parameter family of partial differential equations derived in
\cite{popkovschuetz} essentially differs form our partial differential
equations
(\ref{eq:brickl}). (Actually there is no parameter value for which the
two pde-s are equivalent.)
However, the two families of models  show
many similarities and do have common generalizations.

\section{Sketch of proof}
\label{section:proof}

The present section is divided into four subsections. In subsection
\ref{subs:first} we present the first steps of the `relative entropy
method' applied. As there is no real novelty in this part, we only
list the main steps \emph{without the computational details} which are
essentially the same as in the original work \cite{yau} of Yau or in
Chapter 6. of \cite{kipnislandim} or in \cite{fritz}.

It turns out that for a general two
(or more) component system some identity relating the macroscopic
fluxes $\Phi$ and $\Psi$ is essentially needed for completing the
proof. These relations
are reminiscent of Onsager's reciprocity  relations of nonequilibrium
thermodynamics, see e.g. Chapter 10.D of \cite{reichl}.
However an essential difference is worth noting: while the traditional
Onsager relations are derived under the condition of reversibility of
the microscopic dynamics, in our case conditions
(\ref{cond:cons}), (\ref{cond:staci1}) and (\ref{cond:staci2}) are
involved which do not imply reversability by any means.
Seemingly, these relations were not explicitly noted so far in
the context   of mathematically rigorous hydrodynamic limits. This
omission is probably
explained by the fact that in the concrete models investigated so far
these identities just droped out without any computations.

In subsection
\ref{subs:OR} we prove that
under the conditions (\ref{cond:cons}), (\ref{cond:staci1}) and
(\ref{cond:staci2}) these identities hold in general.
We shall refer to these identities as \emph{Onsager relations}.
It also follows
from these identities that the systems of partial differential
equations arising as hydrodynamic limits under Eulerian scaling are
indeed of hyperbolic type and the thermodynamic equilibrium entropy
$S(u,v)$ is globally convex Lax entropy of the hydrodynamic equations,
as it is commonly assumed. In subsection
\ref{subs:orconseq} we formulate the consequences of the
Onsager relations which are crucial for the further steps of
the proof of the hydrodynamic limit.

Finally, in subsection
\ref{subs:last} we sketch the last steps of the proof. Here again we
follow the standard steps of the relative entropy method, so we omit
all computational details, referring only to the main stations  of
the proof.   For the computational details of subsections
\ref{subs:first} and \ref{subs:last} we refer the reader to Chapter
6. of \cite{kipnislandim} or to \cite{fritz}. However, we warn the
reader that the
omitted details (in particular the last two steps: the control of the
block replacement and the one-block estimate) are rather sophisticated
and mathematically deep.

\subsection{First transformations}
\label{subs:first}

In order to obtain the main estimate (\ref{eq:main})
our aim is to get a Gr\"{o}mwall type inequality:
we will prove that for every $t\in[0,T]$
\beq
\label{eq:gromwall}
\HN(t)- \HN(0) \leq C \int_0^t \HN(s)ds +
o(N),
\eeq
where the error term is uniform in $t\in[0,T]$. Because it is
assumed that $\HN(0)=o(N)$, the Theorem follows.

For proving (\ref{eq:gromwall}) we try to bound (from above)
$\pt \HN(s)$ by $\text{const} \cdot \HN(s)+o(N)$, uniformly for
$s\in[0,T]$. We start form the inequality (\ref{eq:entrprod}) which is
derived in \cite{kipnislandim} under very general conditions, valid in
our case.
\beq
\label{eq:entrprod}
\pt \HN(t)
\leq
N\int_{\OmN}
\frac{\LNa\fN{t}}{ \fN{t}}
d\muN{t}
-
\int_{\OmN}
\frac{\pt \fN{t}}{\fN{t}}
d\muN{t}.
\eeq
Next we transform the two terms appearing on the right hand side
of (\ref{eq:entrprod}). Equation (\ref{eq:l*f/f}) follows form the
smoothness of the functions $\theta(t,x)$ and $\tau(t,x)$ and from
the entropy inequality applied to the measures $\muN{t}$ compared
with the absolute reference measure $\piN{}$.
\beq
\label{eq:l*f/f}
N\int_{\OmN}
\frac{\LNa\fN{t}}{\fN{t}}
d\muN{t}
&=&
-
\sum_{j\in \TN}
\px\theta(t,j/N)
\int_{\OmN} \phi_j d\muN{t}
\\
\nonumber
&&
-
\sum_{j\in \TN}
\px \tau(t,j/N)
\int_{\OmN} \psi_j d\muN{t}
\\
\nonumber
&&
+\Ordo\left(1\right)
\eeq
Equation (\ref{eq:ptf/f}) follows from direct computation of the time
derivative of the function $\fN{t}$.
\beq
\label{eq:ptf/f}
\int_{\OmN}
\frac{\pt\fN{t}}{\fN{t}}
d\muN{t}
&=&
\phantom{+}
\sum_{j\in \TN}
\pt\theta(t,j/N)
\int_{\OmN}
(\xi_j-u(t,j/N))
d\muN{t}
\\
\nonumber
&&
+
\sum_{j\in \TN}
\pt\tau(t,j/N)
\int_{\OmN}
(\eta_j-v(t,j/N))
d\muN{t}
\eeq
Next we replace the local variables $\phi_j$ and $\psi_j$ in
(\ref{eq:l*f/f}), respectively, $\xi_j$ and $\eta_j$ in
(\ref{eq:ptf/f}) by their block averages defined as follows:
if $\delta_j$ is a local microscopic variable its block average is
defined as
\[
\delta_j^l:=\frac{\delta_j+\dots+\delta_{j+l-1}}{l}.
\]
In the following two block-replacements we use again the smoothness of
the functions $\theta(t,x)$ and $\tau(t,x)$ and  the
entropy inequality applied to the measures $\muN{t}$ compared with the
absolute reference measure $\piN{}$.
\beq
\label{eq:l*f/fbl}
N\int_{\OmN}
\frac{\LNa\fN{t}}{\fN{t}}
d\muN{t}
&=&
-
\sum_{j\in \TN}
\px\theta(t,j/n)
\int_{\OmN} \phi^l_j d\muN{t}
\\
\nonumber
&&
-
\sum_{j\in \TN}
\px \tau(t,j/n)
\int_{\OmN} \psi^l_j d\muN{t}
\\
&&
\nonumber
+
\Ordo\left(l\right)
\eeq
\beq
\label{eq:ptf/fbl}
\int_{\OmN}
\frac{\pt\fN{t}}{\fN{t}}
d\muN{t}
&=&
\phantom{+}
\sum_{j\in \TN}
\pt\theta(t,j/N)
\int_{\OmN}
(\xi^l_j-u(t,j/N))
d\muN{t}
\\
\nonumber
&&
+
\sum_{j\in \TN}
\pt\tau(t,j/N)
\int_{\OmN}
(\eta^l_j-v(t,j/N))
d\muN{t}
\\
\nonumber
&&
+
\Ordo\left(l\right)
\eeq
The last transformation of this first, preparatory part is replacing
in (\ref{eq:l*f/fbl}) the block averages $\phi_j^l$, respectively,
$\psi_j^l$ by their equilibrium averages computed at the empirical
densities,
$\Phi(\xi^l_j,\eta^l_j)$, respectively, $\Psi(\xi^l_j,\eta^l_j)$. The
error terms appearing in the third and fourth lines of the right hand
side of (\ref{eq:l*f/fblrepl}) are the most important error  terms to
be controlled by the so called \emph{one block estimate} towards the
end of the proof.
\beq
\label{eq:l*f/fblrepl}
N\int_{\OmN}
\frac{\LNa\fN{t}}{\fN{t}}
d\muN{t}
&=&
-
\sum_{j\in \TN}
\px\theta(t,j/n)
\int_{\OmN} \Phi(\xi^l_j,\eta^l_j) d\muN{t}
\\
\nonumber
&&
-
\sum_{j\in \TN}
\px \tau(t,j/n)
\int_{\OmN} \Psi(\xi^l_j,\eta^l_j) d\muN{t}
\\
\nonumber
&&
-
\sum_{j\in \TN}
\px\theta(t,j/n)
\int_{\OmN}
\left(\phi^l_j-\Phi(\xi^l_j,\eta^l_j)\right) d\muN{t}
\\
\nonumber
&&
-
\sum_{j\in \TN}
\px\tau(t,j/n)
\int_{\OmN}
\left(\psi^l_j-\Psi(\xi^l_j,\eta^l_j)\right) d\muN{t}
\\
\nonumber
&&
+
\Ordo\left(l\right)
\eeq

Before going on with the standard steps of the relative entropy proof
we need to make a detour.

\subsection{An Onsager type identity}
\label{subs:OR}

\begin{lemma}
\label{lem:OR}
Suppose we have a particle system with two conserved quantities
and rates satisfying conditions (\ref{cond:cons}),
(\ref{cond:staci1}) and (\ref{cond:staci2}). Then there exists a
potential function $(\theta,\tau)\mapsto U(\theta, \tau)$ such that
\beq
\label{eq:OR}
\begin{array}{l}
\Phi(\theta,\tau):=
\Phi(u(\theta,\tau),v(\theta,\tau))=
U'_\theta,
\\[8pt]
\Psi(\theta,\tau):=
\Psi(u(\theta,\tau),v(\theta,\tau))=
U'_\tau,
\end{array}
\eeq
or, equivalently
\beq
\label{eq:OR2}
\Phi'_\tau
=
\Psi'_\theta.
\eeq
\end{lemma}

\begin{proof}
We prove (\ref{eq:OR2}).
Throughout the forthcoming proof we adopt the notation
$\xi_j:=\xi(\omega_j)$, $\xi'_j:=\xi(\omega'_j)$, etc.

From the definitions
\[
\ptt(\omega_1)\ptt(\omega_2)=
\exp\{
\theta(\xi_1+\xi_2)
+
\tau(\eta_1+\eta_2)
-
2G(\ttt)
\}
\pi(\omega_1)\pi(\omega_2),
\]
and
\beqs
&&
{\big(\ptt(\omega_1)\ptt(\omega_2)\big)}'_\theta
=
\\[8pt]
&&
\hskip10mm
\pi(\omega_1)\pi(\omega_2)
e^{\txte-2G(\ttt)}
\left\{(\xi_1+\xi_2)-2 u(\ttt) \right\}
=
\\[8pt]
&&
\hskip10mm
\frac{\pi(\omega_1)\pi(\omega_2)}{Z(\ttt)^{3}}
\sum_{\omega_3\in S}
\pi(\omega_3)
(\xi_1+\xi_2-2\xi_3)
e^{
\theta(\xi_1+\xi_2+\xi_3)
+
\tau(\eta_1+\eta_2+\eta_3)
},
\eeqs
where $Z(\ttt)=\exp\{G(\ttt)\}$.  Similarly,
\beqs
&&
{\big(\ptt(\omega_1)\ptt(\omega_2)\big)}'_\tau
=
\\[8pt]
&&
\hskip10mm
\frac{\pi(\omega_1)\pi(\omega_2)}{Z(\ttt)^{3}}
\sum_{\omega_3\in S}
\pi(\omega_3)
(\eta_1+\eta_2-2 \eta_3)
e^{
\theta(\xi_1+\xi_2+\xi_3)
+
\tau(\eta_1+\eta_2+\eta_3)
}.
\eeqs
Hence
\beqs
\Phi'_\tau(\ttt)
&=&
\frac{1}{Z(\ttt)^3}
\sum_{\omega_1,\omega_2,\omega_3\in S}
\pi(\omega_1)\pi(\omega_2)\pi(\omega_3)
e^{
\theta(\xi_1+\xi_2+\xi_3)
+
\tau(\eta_1+\eta_2+\eta_3)
}
\\
&&
\phantom{MMM}
\times\,
(\eta_1+\eta_2-2\eta_3)
\sum_{\omega'_1,\omega'_2\in S}
r(\omega_1,\omega_2,\omega'_1,\omega'_2)
(\xi'_2-\xi_2),
\\
\Psi'_\theta(\ttt)
&=&
\frac{1}{Z(\ttt)^3}
\sum_{\omega_1,\omega_2,\omega_3\in S}
\pi(\omega_1)\pi(\omega_2)\pi(\omega_3)
e^{
\theta(\xi_1+\xi_2+\xi_3)
+
\tau(\eta_1+\eta_2+\eta_3)
}
\\
&&
\phantom{MMM}
\times\,
(\xi_1+\xi_2-2\xi_3)
\sum_{\omega'_1,\omega'_2\in S}
r(\omega_1,\omega_2,\omega'_1,\omega'_2)
(\eta'_2-\eta_2),
\eeqs
For the proof of the lemma it is enough to prove for any
$K\in[3\min\xi,3\max\xi]$ and $L\in[3\min\eta,3\max\eta]$
\beqs
&&
\hskip-8mm
\sum_{\substack{\omega_1,\omega_2,\omega_3,\omega'_1,\omega'_2 \in S:\\
                \xi_1+\xi_2+\xi'=K\\
                \eta_1+\eta_2+\eta'=L
               }
     }
\pxya
r(\omega_1,\omega_2;\omega'_1,\omega'_2)
\\[-20pt]
&&
\hskip22mm
\times\,
\big(
(\eta_1+\eta_2-2\eta_3)
(\xi'_2-\xi_2)
-
(\xi_1+\xi_2-2\xi_3)
(\eta'_2-\eta_2)
\big)
=0
\eeqs
From condition (\ref{cond:cons}) imposed on the rate functions it
follows that in all nonzero terms of the above sum one can replace
$\eta_1+\eta_2$ by $\eta'_1+\eta'_2$ and  $\eta'_2-\eta_2$ by
$\eta_1-\eta'_1$,  and similarly for the $\xi$-s. After
straightforward computations this equation  (still to be proved)
becomes
\beqs
&&
\hskip-20mm
\sum_{\substack{\omega_1,\omega_2,\omega_3,\omega'_1,\omega'_2 \in S:\\
                \xi_1+\xi_2+\xi'=K\\
                \eta_1+\eta_2+\eta'=L
               }
     }
\pxya
r(\omega_1,\omega_2;\omega'_1,\omega'_2)
\\[-20pt]
&&
\hskip20mm
\times\,
\big(
\Delta(\omega_1,\omega_2,\omega_3)
+
\Delta(\omega'_2,\omega'_1,\omega_3)
\big)
=0
\eeqs
where $\Delta:S\times S\times S \to Z$ is defined as follows
\beqs
\Delta(\omega_1,\omega_2,\omega_3)
:=
\xi_1(\eta_2-\eta_3)
+
\xi_2(\eta_3-\eta_1)
+
\xi_3(\eta_1-\eta_2).
\eeqs
Note that $\Delta$ is \emph{antisymmetric} regarding permutation of
its variables.

Next, from condition (\ref{cond:staci1}) it follows that the
contribution of the two terms on the left hand side of the
previous equation is the same. Thus, it is to be proved that
\beqs
&&
\hskip-12mm
\sum_{\substack{\omega_1,\omega_2,\omega_3,\omega'_1,\omega'_2 \in S:\\
                \xi_1+\xi_2+\xi'=K\\
                \eta_1+\eta_2+\eta'=L
               }
     }
\pxya
r(\omega_1,\omega_2;\omega'_1,\omega'_2)
\Delta(\omega_1,\omega_2,\omega_3)
=
\\[-20pt]
&&
\hskip20mm
\sum_{\substack{\omega_1,\omega_2,\omega_3 \in S:\\
                \xi_1+\xi_2+\xi'=K\\
                \eta_1+\eta_2+\eta'=L
               }
     }
\pxya
R(\omega_1,\omega_2)
\Delta(\omega_1,\omega_2,\omega_3)
=
0
\eeqs

Finally, from the antisymmetry of the function $\Delta$ and condition
(\ref{cond:staci2}) imposed on the function $R$ it follows indeed that
this last sum equals zero.

\end{proof}

\subsection{Consequences of the Onsager relations}
\label{subs:orconseq}

Relation (\ref{eq:OR2})  is the same as saying that
the matrix $D(u(\theta,\tau),v(\theta,\tau))\cdot G''(\theta,\tau)$ is
symmetric. Using (\ref{eq:inverse}) this also reads as
\beq
\label{eq:sym}
S''(u,v)\cdot  D(u,v)
=
\left(
S''(u,v)\cdot D(u,v)
\right)^{\dagger}.
\eeq
This relation implies that only \emph{hyperbolic} two-by-two
systems of conservation laws (\ref{eq:pde}) can arise as
hydrodynamic limits. Indeed, as the following elementary
argument  shows relation (\ref{eq:sym}) can hold with a positive
definite matrix $S''$ only if $D(u,v)$ can be diagonalized (in the
real sense),
which is exactly the condition of hyperbolicity of the system
(\ref{eq:pde}).
Indeed, since $S''$ is positive definite, we can write
%
%
%
%
\beq \label{eq:diag} D=(S'')^{-1/2} \left((S'')^{-1/2} (S''D)
(S'')^{-1/2}\right) (S'')^{1/2}, \eeq which means that $D$ is
similar to  the real  symmetric matrix 
 $(S'')^{-1/2}(S''D)
(S'')^{-1/2}$, and from this the (real) diagonalizability of $D$
follows.
Furthermore, (\ref{eq:sym})  is spelled out as
\beq
\label{eq:laxentropy}
S''_{uu}\Phi'_v  + S''_{uv}\Psi'_v
=
S''_{vu}\Phi'_u  + S''_{vv}\Psi'_u,
\eeq
which is readily recognized as the the partial differential
equation defining the \emph{Lax entropies} of the system
(\ref{eq:pde}). The function $F(u,v):=U(\theta(u,v),\tau(u,v))$ is
the corresponding (macroscopic) entropy-flux.  See \cite{serre} or
\cite{smoller} for the pde notions involved. Thus, part (i) of the
Theorem is proved.

Now we turn to two further consequences of Lemma \ref{lem:OR}
which turn out to be of crucial importance in the hydrodynamic
behavior.

First, the time derivatives of $\theta$ and $\tau$ are expressed.
From the pde (\ref{eq:pde}) it follows that
\beqs
\pt \theta
=
-\theta'_u \Phi'_u \px u
-\theta'_u \Phi'_v \px v
-\theta'_v \Psi'_u \px u
-\theta'_v \Psi'_v \px v.
\eeqs
Using the identity (\ref{eq:laxentropy}) we replace
\beqs
\theta'_u \Phi'_v
=
\theta'_v\Phi'_u
+
\tau'_v\Psi'_u
-
\tau'_u\Psi'_v
\eeqs
in the second term of the right hand side.
Using also the straightforward identities
$u'_\tau=v'_\theta$ and $\theta'_v=\tau'_u$
(see subsection \ref{subs:expectations})
we finally get
\beq
\label{eq:pttheta}
\pt \theta
=
\Phi'_u\px\theta + \Psi'_u\px\tau,
\eeq
and by identical considerations
\beq
\label{eq:pttau}
\pt \tau
=
\Phi'_v\px\theta + \Psi'_v\px\tau.
\eeq
Second, due to identity (\ref{eq:OR}),
\beq
\nonumber
&
\dsum_{j\in\TN}
\Big(
\px\theta(j/N)\Phi(u(j/N),v(j/N))
+
\px\tau(j/N)\Psi(u(j/N),v(j/N))
\Big)
\\
&
\qquad\qquad
=
\label{eq:complete}
\dsum_{j\in\TN}
\px U(u(j/N),v(j/N))
=
\Ordo(1).
\eeq

\subsection{End of proof}
\label{subs:last}

Now we return to proving (\ref{eq:gromwall}).
Denote
\beqs
{\cal D}\Phi(u,v;\tilde u,\tilde v)
:=
\Phi(\tilde u,\tilde v) - \Phi(u,v)
-\Phi'_u(u,v)(\tilde u-u)
-\Phi'_v(u,v)(\tilde v-v)
\eeqs
and similarly for ${\cal D}\Psi(u,v;\tilde u,\tilde v)$.
Applying (\ref{eq:pttheta}), (\ref{eq:pttau}) and (\ref{eq:complete}),
from (\ref{eq:ptf/fbl}) and (\ref{eq:l*f/fblrepl}) we obtain
\beq
\label{eq:sumup}
&&
\hskip-1cm
\int_{\OmN}
\frac{\pt \fN{t}-N\LNa\fN{t}}{ \fN{t}}
d\muN{t}
=
\\
\nonumber
&&
\hskip1cm
\phantom{+}
\sum_{j\in \TN}
\px\theta(t,j/n)
\int_{\OmN}
{\cal D}\Phi(u(t,j/N),v(t,j/N);\xi^l_j,\eta^l_j)
d\muN{t}
\\
\nonumber
&&
\hskip1cm
+
\sum_{j\in \TN}
\px \tau(t,j/n)
\int_{\OmN}
{\cal D}\Psi(u(t,j/N),v(t,j/N);\xi^l_j,\eta^l_j)
d\muN{t}
\\
\nonumber
&&
\hskip1cm
+
\sum_{j\in \TN}
\px\theta(t,j/n)
\int_{\OmN}
\left(\phi^l_j-\Phi(\xi^l_j,\eta^l_j)\right) d\muN{t}
\\
\nonumber
&&
\hskip1cm
+
\sum_{j\in \TN}
\px\tau(t,j/n)
\int_{\OmN}
\left(\psi^l_j-\Psi(\xi^l_j,\eta^l_j)\right) d\muN{t}
\\
\nonumber && \hskip1cm + \Ordo\left(l\right) \eeq
The first two terms on the right hand side of (\ref{eq:sumup}) are
estimated by the
entropy inequality, comparing the measure $\muN{t}$ with the
\emph{local equilibrium measure} $\nuN{t}$:
\beq
\label{eq:replace}
&&
\hskip-12mm
\sum_{j\in \TN}
\int_{\OmN}
\left(\left|
{\cal D}\Phi(u(t,j/N),v(t,j/N);\xi^l_j,\eta^l_j) \right|\right.
\\
\nonumber
&&
\hskip10mm
\left.
+
\left|{\cal D}\Psi(u(t,j/N),v(t,j/N);\xi^l_j,\eta^l_j)
\right|\right)
d\muN{t}
\\
\notag
&&
\hskip65mm
\le
C\HN(t) +\Ordo(Nl^{-1}).
\eeq
The last two terms in (\ref{eq:sumup}) are estimated only
\emph{integrated against time}. Applying the so-called one block
estimate (see e.g. Chapter 5 of \cite{kipnislandim}) one gets
\beq
\nonumber
&&
\lim_{l\to\infty}\lim_{N\to\infty}
N^{-1}\sum_{j\in \TN}
\int_0^tds
\int_{\OmN}
\abs{\phi^l_j-\Phi(\xi^l_j,\eta^l_j)} d\muN{t}=
0,
\\
\label{eq:obe}
\\
\nonumber
&&
\lim_{l\to\infty}\lim_{N\to\infty}
N^{-1}\sum_{j\in \TN}
\int_0^tds
\int_{\OmN}
\abs{\psi^l_j-\Psi(\xi^l_j,\eta^l_j)} d\muN{t}=
0.
\eeq
This is the only part of the proof where condition
(\ref{cond:irred}) is used, which ensures ergodicity of the Markov
process $\XN_t$ on the `hyperplanes' $\OmN_{K,L}$.

Finally, inserting  (\ref{eq:replace}) and  (\ref{eq:obe}) in
(\ref{eq:sumup}), via (\ref{eq:entrprod}) we obtain
(\ref{eq:gromwall}) and  thus the part (ii) of the Theorem is also
proved.

\section{Particle systems with several conserved variables}
\label{section:several}

As noted in the introduction, the results described in the
previous sections are also valid for particle systems with more
than 2 conserved quantities.

Before we formulate the general results we have to summarize some
notations . Let $n\geq 2$ be fixed integer, and
$\vct{\xi}=(\xxi{1},\xxi{2},\dots,\xxi{n}): S \rightarrow \R^n$
the vector of conserved quantities. Throughout the present section
bold face symbols will denote $n$-vectors.

We require the rate function
to satisfy similar conditions as listed in subsection
\ref{subs:rates} (in place of conditions (\ref{cond:cons}) and
(\ref{cond:irred}) we need the suitable generalizations). For
every $\vct{\theta}\in \R^n$ we can define momentum generating
function $G(\vct{\theta})$ as
\[
G(\vct{\theta})
:=
\log \sum_{\omega \in S}
e^{\vct{\theta}\cdot\vct{\xi}(\omega)} \pi(\omega),
\]
and the probability measures
\beqs
\pi_{\vct{\theta}}(\omega)
&:=&
\pi(\omega)
\exp(\vct{\theta}\cdot\vct{\xi}(\omega)-G(\vct{\theta}))
\\
\piN{\vct{\theta}}
&:=&
\prod_{j\in \TN} \pi_{\vct{\theta}}
\eeqs
on $S$, respectively, on $\OmN$. We define the
expectation of the conserved quantities with respect to the
measure $\piN{\vct{\theta}}$:
\[
\vct{u}(\vct{\theta})
:=
\expect_{\vct{\theta}}(\vct{\xi})
=
\nabla_{_{\!\!\vct{\theta}}} G(\vct{\theta}).
\]
One can easily show,  that
$\nabla^2_{_{\!\!\vct{\theta}}} G(\vct{\theta})
= \cov_{\vct{\theta}}(\vct{u},\vct{u})$
is positive definite. As a consequence,
the function $\vct{\theta}\mapsto \vct{u}$
is invertible,  and
\[
\vct{\theta}(\vct{u})
=
\nabla_{_{\!\!\vct{u}}} S(\vct{u}),
\]
where $S(\vct{u})$ is the convex conjugate of $G(\vtt)$:
\[
S(\vtu)
:=
\sup_{\vtu \in \R^n} (\vtu \cdot \vtt -G(\vtt)).
\]
We introduce the flux of the vector of the conserved quantities
and its expectation:
\[
\begin{array}{lcl}
\vct{\phi}(\omega_1,\omega_2)
&:=&
\dsum_{\omega'_1,\omega'_2\in S}
r(\omega_1,\omega_2;\omega'_1,\omega'_2)
(\vct{\xi}(\omega'_2)-\vct{\xi}(\omega_2))
\\[10pt]
\vct{\Phi}(\vtu)&
:=&
\expect_{\vtt(\vtu)} \vct{\phi}
\\[10pt]
&=&
\dsum_{\substack{\omega_1,\omega_2,\\ \omega'_1,\omega'_2}\in S}
r(\omega_1,\omega_2;\omega'_1,\omega'_2)
(\vct{\xi}(\omega'_2)-\vct{\xi}(\omega_2))
\pi_{\vtu}(\omega_1)\pi_{\vtu}(\omega_2).
\end{array}
\]
Now we are able to formulate  results of the previous sections in
the more general setting.

Using the arguments presented in section \ref{section:proof} one
can show that under Eulerian scaling the vector of the local
densities of the conserved quantities $\vtu(t,x)$ evolve according
to the following $n$-component partial differential equation:
\beq
\label{eq:pde_gen}
\pt \vtu+\px \vct{\Phi} (\vtu)=0.
\eeq

Lemma \ref{lem:OR} applies for any two conserved quantities
$\xxi{i}, \xxi{j}$ ($i \neq j$), thus if we denote the derivative
matrix of the flux vector $\vct{\Phi}(\vtu)$ by
$D(\vtu):=\nabla_{_{\!\!\vct{u}}}\vct{\Phi}(\vtu)$ and the second derivative
matrix of
the thermodynamic entropy
$S''(\vct{u}):=\nabla^2_{_{\!\!\vct{u}}}S(\vct{u})$, we get
\beq
\label{eq:sdsym}
S''(\vtu) \cdot D(\vtu)
=
\left( S''(\vtu) \cdot D(\vtu)\right)^\dagger
\eeq
Since $S''(\vtu)$ is positive
definite (\ref{eq:diag}) implies that $D(\vtu)$ can be
diagonalized which
means that the arising system of partial differential equations is
\emph{hyperbolic}. Moreover, (\ref{eq:sdsym}) spelled out is
\beq
\label{eq:sevlaxentropy}
\frac{\partial^2S}{\partial u_i \partial u_i}
\frac{\Phi_i}{\partial u_j}
+
\frac{\partial^2S}{\partial u_i \partial u_j}
\frac{\Phi_j}{\partial u_j}
=
\frac{\partial^2S}{\partial u_j \partial u_i}
\frac{\Phi_i}{\partial u_i}
+
\frac{\partial^2S}{\partial u_j \partial u_j}
\frac{\Phi_j}{\partial u_i},
\eeq
with $1\le i<j\le n$.
These are exactly the $n(n-1)/2$ equations defining the Lax entropies
of the hyperbolic system (\ref{eq:pde_gen}).
It is well-known that in the case of $n\geq 3$ only very special
$n$-component hyperbolic conservation laws possess Lax entropies.
In general, the defining equations (\ref{eq:sevlaxentropy})
are overdetermined. In \cite{serre} these paricular systems
of hyperbolic coservation laws  are
called of \emph{physical} type. From the previous arguments
it follows that only physical hyperbolic equations can arise as the
hydrodynamic limit of an interacting particle system satisfying our
conditions.


\bigskip

\noindent {\bf\large Acknowledgement:} The authors wish to thank
J\'ozsef Fritz for the many illuminating discussions on the topics
of hydrodynamic limits and for pointing out Corollary (ii). We also
thank J\'anos Kert\'esz for guiding comments on Onsager's relations.


\vskip1cm

\hbox{\sc
\vbox{\noindent
\hsize66mm
B\'alint T\'oth\\
Institute of Mathematics\\
Technical University Budapest\\
Egry J\'ozsef u. 1.\\
H-1111 Budapest, Hungary\\
{\tt balint{@}math.bme.hu}
}
\hskip5mm
\vbox{\noindent
\hsize66mm
Benedek Valk\'o\\
Institute of Mathematics\\
Technical University Budapest\\
Egry J\'ozsef u. 1.\\
H-1111 Budapest, Hungary\\
{\tt valko{@}math.bme.hu}
}
}


\begin{thebibliography}{99}


\bibitem{balazs}
M. Bal\'azs: Growth fluctuations in interface models. {\sl
Preprint} (2001)

\bibitem{cocozza}
C. Cocozza: Processus des misanthropes. {\sl Zeitschrift f\"ur
Wahrscheinlichkeitstheorie und verwandte Gebiete} {\bf 70}:
509-523 (1985)

\bibitem{fritz}
J. Fritz: {\sl An Introduction to the Theory of Hydrodynamic
Limits.\/} Lectures in Mathematical Sciences. Graduate School of
Mathematics, Univ. Tokyo, 2001.


\bibitem{kipnislandim}
C. Kipnis, C. Landim: {\sl Scaling Limits of Interacting Particle
Systems.\/} Springer, 1999.

\bibitem{lax}
P. Lax: Shock waves and entropy. In: {\sl Contributions to
  Nonlinear Functional Analysis}, ed.: E.A. Zarantonello. 
Academic Press, 1971, pp. 603-634

\bibitem{popkovschuetz}
V. Popkov, G.M. Sch\"utz:
Shocks and excitation dynamics in driven diffusive two channel
systems. {\sl Preprint} (2002)

\bibitem{reichl}
L.E. Reichl:
{\sl A Modern Course in Statistical Physics}, Second Edition,
John Wiley and Sons, 1998

\bibitem{rezakhanlou2}
F. Rezakhanlou: Microscopic structure of shocks in one
conservation laws. {\sl Annales de l'Institut Henri Poincar\'e ---
Analyse Non Lineaire} {\bf 12}: 119-153 (1995)

\bibitem{serre}
D. Serre: {\sl Systems of Conservation Laws.} Vol 1-2.
Cambridge University Press, 2000

\bibitem{smoller}
J. Smoller: {\sl Shock Waves and Reaction Diffusion Equations}, Second
Edition, Springer, 1994.

\bibitem{tothwerner}
B. Tóth, W. Werner: Hydrodynamic equation for a deposition model,
en collaboration, In: {\sl In and out of equilibrium. Probability
with a physics flavor,} V. Sidoravicius Ed., Progress in
Probability {\bf 51}, Birkhäuser, 227-248 (2002)

\bibitem{tothvalko}
B. T\'oth, B. Valk\'o:
Between equilibrium fluctuations and Eulerian scaling. Perturbation of
equilibrium for a class of deposition models. {\sl Journal of
  Statistical Physics} {\bf 109}: 177-205 (2002)

\bibitem{yau}
H.T. Yau: Relative entropy and hydrodynamics of Ginzburg-Landau
models. {\sl Lett. Math. Phys.\/} {\bf 22}: 63-80 (1991)

\end{thebibliography}
\end{document}